\documentclass[oneside,draft,11pt]{amsart}

\newcommand{\R}{\ensuremath{\mathbb{R}}}
\newcommand{\dif}{{\rm d}}
\newcommand{\be}{\begin{equation}}
\newcommand{\ee}{\end{equation}}
\newcommand{\beq}{\begin{eqnarray}}
\newcommand{\eeq}{\end{eqnarray}}
\newcommand{\wt}{\widetilde}
\newcommand{\ol}{\overline}

\newcommand{\Id}{\mathrm{Id}}
\newcommand{\g}{{\mathfrak{g}}}
\newcommand{\h}{{\mathfrak{h}}}
\newcommand{\A}{{\mathfrak{a}}}
\newcommand{\m}{{\mathfrak{m}}}
\newcommand{\te}{{\mathfrak{t}}}
\newcommand{\gls}{\mathfrak{gl}}
\newcommand{\gl}{{\mathfrak{gl}(n)}}
\newcommand{\gx}{{\mathfrak{g}_x}}
\newcommand{\glTx}{{\mathfrak{gl}(T_xM)}}
\newcommand{\tripla}{{(M,\nabla,P)}}
\newcommand{\bfd}{{\mathrm{d\!I}}}
\newcommand{\Dd}{\mathrm D}
\newcommand{\Ad}{\mathrm{Ad}}
\newcommand{\ad}{\mathrm{ad}}
\newcommand{\I}{{{\mathfrak{I}}}}
\newcommand{\FR}{{{\rm{FR}}}}
\newcommand{\GL}{{{\rm{GL}}}}
\newcommand{\Lin}{\mathrm{Lin}}
\newcommand{\Alt}{\mathrm{Alt}}
\newcommand{\pro}[1]{{\mathfrak p_{\mathrm{#1}}}}

\newcommand{\bfG}{\mathbf{\Gamma}}
\newcommand{\pip}{{\Pi}}
\newcommand{\ttimes}{\mathbin{\hbox to 0pt{$\times$\hskip -6.7pt\lower2.8pt\hbox{$\scriptscriptstyle\sim$}\hskip 0pt minus 1fil}\hphantom\times}}
\newcommand{\Cat}[1]{\underline{\mathfrak{#1}}}
\newcommand{\Vect}[1]{\underline{\mathfrak{Vec}}^{#1}}

\usepackage[english]{babel}
\usepackage[latin1]{inputenc}
\usepackage{times}
\usepackage{dsfont}
\usepackage[all]{xy}

\title[An algebraic characterization of affine manifolds with $\mathbf G$-stucture]%
{An algebraic characterization of affine manifolds with $\mathbf G$-stucture satisfying a homogeneity condition}

\author[C. Marín]{Carlos Alberto Marin Arango}

\address{Departamento de Matem\'atica,\hfill\break\indent  Universidade de Antioquia, Colombia}
\email{camara@matematicas.udea.edu.co}



\date{August 2010}

\begin{document}

\makeatletter
\renewenvironment{proof}[1][\proofname]{\par
  \pushQED{\qed}%
  \normalfont \topsep6\p@\@plus6\p@\relax
  \trivlist
  \item[\hskip\labelsep
        \scshape
    #1\@addpunct{.}]\ignorespaces
}{%
  \popQED\endtrivlist\@endpefalse
}
\makeatother

\numberwithin{equation}{section}
\theoremstyle{remark}\newtheorem{obs}{Remark}
\theoremstyle{plain}\newtheorem{teo}{Theorem}[section]
\theoremstyle{plain}\newtheorem{prop}[teo]{Proposition}
\theoremstyle{plain}\newtheorem{lema}[teo]{Lemma}
\theoremstyle{plain}\newtheorem{cor}[teo]{Corollary}
\theoremstyle{definition}\newtheorem{df}[teo]{Definition}
\theoremstyle{remark}\newtheorem{rem}[teo]{Remark}
\swapnumbers
\theoremstyle{definition}\newtheorem{example}[teo]{Example}

\begin{abstract}
We give an algebraic characterization of the possible characteristic tensors of an infinitesimally
homogeneous affine manifold with $G$-structure. Such concepts were introduced in \cite{PiTausk}.
\end{abstract}

\maketitle

\begin{section}{Introduction}

The concept of {\em infinitesimally homogeneous\/} affine manifold with $G$-structure was introduced in
the recent article \cite{PiTausk} with the aim to find a unifying language for several isometric
immersion (Bonnet type) theorems that appear in the classical literature \cite{Dajczer}
(immersions into Riemannian manifolds with constant sectional curvature,
immersions into K\"ahler manifolds of constant holomorphic curvature),
and also some more recent results (see for instance \cite{Benoit1, Benoit}) concerning
the existence of isometric immersions in more general Riemannian
manifolds. By an affine manifold with $G$-structure we mean a triple $(M,\nabla,P)$, with $M$ an $n$-dimensional
diffe\-rentiable manifold, $\nabla$ a connection on $M$ and $P$ a $G$-structure on $M$, i.e., $G$ is a Lie subgroup
of $\GL(n)$ and $P$ is a $G$-principal subbundle of the frame bundle of $M$. We denote by $R$ and $T$, respectively, the
curvature and torsion tensors of $\nabla$. In order to handle the case in which $P$ is not compatible with $\nabla$,
the concept of {\em inner torsion\/} was introduced in \cite{PiTausk}: it is a tensor $\I^P$ that plays the role of
a covariant derivative of the $G$-structure $P$ and it vanishes if and only if $\nabla$ is compatible with $P$.
The concept of infinitesimal homogeneity plays the same role in the theory of affine manifolds with
$G$-structure as the concept of constant sectional curvature plays in Riemannian geometry; in fact, Riemannian
manifolds with constant sectional curvature are precisely the infinitesimally homogeneous triples $(M,\nabla,P)$
in which $P$ is the $\mathrm O(n)$-principal bundle of orthonormal frames and both the torsion and the inner torsion
vanish. Notice that Riemannian manifolds with constant sectional curvature are those in which the (four indexed) matrix
representing the curvature tensor with respect to orthonormal frames is independent of the orthonormal frame and of
the point on the manifold. While it does not make sense to require that a tensor field on a manifold be constant,
we can define, for manifolds endowed with a $G$-structure, the notion of {\em $G$-constant\/} tensor field: that
is a tensor field whose matrix with respect to frames that belong to the $G$-structure is independent of the
frame and of the point of the manifold. An affine manifold with $G$-structure $(M,\nabla,P)$ is said to be
{\em infinitesimally homogeneous\/} if the tensor fields $R$, $T$ and $\I^P$ are all $G$-constant.
When $M$ is simply connected and $\nabla$ is geodesically complete then this condition implies that the group
of all affine $G$-structure preserving diffeomorphisms of $M$ acts transitively on the frames that belong to $P$
and in that case we say that the triple $(M,\nabla,P)$ is {\em homogeneous\/} \cite{PiTausk}.

The $G$-constant tensor fields $R$ and $T$ of an infinitesimally homogeneous triple $(M,\nabla,P)$ are represented,
with respect to an arbitrary frame belonging to $P$, by multilinear maps $R_0:\R^n\times\R^n\times\R^n\to\R^n$ and
$T_0:\R^n\times\R^n\to\R^n$, respectively; moreover, the $G$-constant inner torsion $\I^P$ is represented
(with respect to an arbitrary frame belonging to $P$) by a linear map $\I^P_0:\R^n\to\gl/\mathfrak g$,
where $\mathfrak g$ denotes the Lie algebra of $G$. We call $R_0$, $T_0$, $\I^P_0$ the {\em characteristic tensors\/}
of $(M,\nabla,P)$. The characteristic tensors $R_0$, $T_0$, $\I^P_0$ characterize locally an
infinitesimally homogeneous triple $(M,\nabla,P)$, in the sense that two infinitesimally homogeneous triples
having the same characteristic tensors are locally equivalent (by means of affine $G$-structure preserving
diffeomorphisms). It is then very natural to ask what are the necessary and sufficient conditions for maps
$R_0$, $T_0$, $\I^P_0$ to be the characteristic tensors of an infinitesimally homogeneous triple $(M,\nabla,P)$.
This paper answers such question.

The main result of this paper can be seen as part of a program of reducing a problem of classification
of certain geometric objects to a problem of classification of certain algebraic objects. Other examples
of such reductions are: (i) the result that two Lie groups having the same Lie algebra are locally isomorphic
and every Lie algebra is the Lie algebra of a Lie group; (ii) the result that two Riemannian symmetric spaces having the
same orthogonal involutive Lie algebra (oil algebra) are locally isometric and every oil algebra is the oil algebra
of a Riemannian symmetric space (see \cite{Helgason}).

It would be natural to ask what are the necessary and sufficient conditions for $R_0$, $T_0$, $\I^P_0$ to be
the characteristic tensors of a (globally) homogeneous triple $(M,\nabla,P)$. Is it true that if $R_0$,
$T_0$, $\I^P_0$ are the characteristic tensors of an infinitesimally homogeneous triple then they are also the
characteristic tensors of some (globally) homogeneous triple? While we do not know the answer to that question,
a partial answer will be given in a forthcoming paper.

\end{section}

\begin{section}{Notation and Preliminaries}

\begin{subsection}{Vector spaces}Let $V$ be a real finite-dimensional vector space. We denote by $\GL(V)$ the general linear group of $V$ and by $\gls(V)$ its
Lie algebra. If $W$ is another real finite-dimensional vector space, $\Lin_k(V;W)$ denotes the space of $k$-linear maps from $V$ to $W$. Given multilinear maps $T\in \Lin_k(V;V)$, $S\in \Lin_k(W;W)$ and a (not necessarily invertible) linear map $\sigma:V\to W$ then $T$ is said to be {\em $\sigma$-related\/} with $S$ if: \[S\big(\sigma (v_1),\dots,\sigma(v_k)\big)=\sigma\big(T(v_1,\dots,v_k)\big),\] for all $v_1,\dots,v_k \in V$. If $p:V\to W$ is a linear isomorphism we denote by $\mathcal I_p:\GL(V)\to\GL(W)$ the Lie group isomorphism given by conjugation with $p$ and $\Ad_p=\dif \mathcal I_p(\Id):\gls(V)\to\gls(W)$ denotes the Lie algebra
isomorphism given by conjugation with $p$.
\end{subsection}

\begin{subsection}{$\mathbf G$-structures on manifolds}

If $G$ is a Lie subgroup of ${\rm GL}(n)$, by a {\em $G$-structure\/} on an $n$-dimensional real vector space $V$ we mean a $G$-orbit of the action given by right composition of ${\rm GL}(n)$ on the set of all linear isomorphisms $p:\R^n\to V$. By a $G$-structure on an $n$-dimensional differentiable manifold $M$ we mean a $G$-principal subbundle $P$ of $ \FR(TM)$, such that for each $x\in M$, $P_x$ is a $G$-structure on the vector space $T_xM$. Let $M$ and $M'$ be $n$-dimensional diffe\-rentiable manifolds endowed with $G$-structures $P$ and $P'$, respectively. A smooth map $f:M\to M'$ is said to be {\em $G$-structure preserving\/} if for each $x\in M$, the linear map $\dif f_x:T_xM\to T_{f(x)}M'$ sends frames of $P_x$ to frames that belong to $P'_{f(x)}$.

\begin{obs}\label{thm:tensorG-invariante}
If $G$ is a Lie subgroup of ${\rm GL}(n)$ a multilinear map $\tau_0 \in \Lin_k(\R^n;\R^n)$ is said to be {\em $G$-invariant\/}, if for each $g \in G$, $\tau_0$ is $g$-related with itself. Clearly, given a $G$-invariant tensor $\tau_0 \in \Lin_k(\R^n;\R^n)$ one can induce a version of $\tau_0$ on every vector space endowed with a $G$-structure. More precisely, let $V$ be a real $n$-dimensional vector space endowed with a $G$-structure $P$. Given any $p\in P$ let $\tau_V\in\Lin_k(V;V)$ be the tensor which is $p$-related with $\tau_0$. The $G$-invariance of $\tau_0$ implies that $\tau_V$ does not depend on the choice
of $p \in P$. In particular, when $M$ is an $n$-dimensional differentiable manifold endowed with a $G$-structure $P$ and $\tau_0 \in \Lin_k(\R^n;\R^n)$
is $G$-invariant, by using frames that belong to $P$ it is possible to define a tensor field $\tau$ on $M$ such that for each $x\in M$, the map $\tau_x\in\Lin_k(T_xM;T_xM)$ is the version of $\tau_0$ in $T_xM$.
\end{obs}
\end{subsection}
\begin{subsection}{Connections on vector bundles}

Let $E$ be a vector bundle over a differentiable manifold $M$ with typical fiber $E_0$. We denote by $\bfG(E)$ the set of all smooth sections of $E$ and by $\FR_{E_0}(E)$ the $\GL(E_0)$-principal bundle over $M$ formed by all $E_0$-frames
of $E$. When $E_0=\R^n$ we write $\FR(E)$ instead of $\FR_{E_0}(E)$. If $\epsilon:U\to E$ is a local section of the vector bundle $E$ and $s:U\to \FR_{E_0}(E)$
is a smooth local frame for $E$ then the {\em representation\/} of the section $\epsilon$ with respect to the smooth local frame $s$ is a map $\tilde{\epsilon}:U\to E_0$ defined by: $\tilde\epsilon(x)=s(x)^{-1}\big(\epsilon(x)\big)$, for
all $x\in U$.

A smooth local frame $s:U\to \FR_{E_0}(E)$ defines, in a natural way, a connection $\bfd^s$ in $E\vert_U$, which corresponds via the trivialization of $E\vert_U$ definided by $s$ to the standard derivative. More explicitly, we set: \[\bfd^s_v\epsilon=s(x)\big(\dif\tilde\epsilon_x(v)\big),\]
for all $x\in U$, $v\in T_xM$ and all $\epsilon\in\bfG(E\vert_U)$, where
$\tilde\epsilon:U\to E_0$ denotes the representation of $\epsilon$ with respect to the local frame $s$.

If $\nabla$ is a connection in $E$, the {\em Christoffel tensor\/} of $\nabla$ with respect to the smooth local frame $s$ is the smooth tensor $\Gamma=\nabla-\bfd^s \in \bfG(TM^*\otimes E^*\otimes E)$ such that:
\[\nabla_v\epsilon=\bfd^s_v\epsilon+\Gamma_x\big(v,\epsilon(x)\big),\]
for all $x\in U$, $v\in T_xM$ and all $\epsilon\in\bfG(E\vert_U)$.
Denoting by $\omega$ the smooth $\gls(E_0)$-valued connection form on $\FR_{E_0}(E)$ associated to $\nabla$, we have the following:
\be\label{eq:tensorcriseforma} \Gamma_x(v)=s(x)\circ\bar\omega_x(v)\circ
s(x)^{-1}\in\gls(E_x), \ee for all $x\in U$, $v\in T_xM$, where $\bar\omega=s^*\omega$ denotes the pullback by $s$ of the connection form $\omega$.

\begin{obs}\label{thm:diferenciadeconexiones}
If $\nabla$ is a (symmetric) connection on $TM$ and
$\te:\bfG(TM)\times\bfG(TM)\to\bfG(TM)$ is an arbitrary $C^\infty(M)$-bilinear (symmetric)
map, $\nabla'=\nabla+\te$ is also a (symmetric) connection on $TM$ and a simple calculation shows that, (see \cite{T}):
\begin{align}
\label{eq:curvaturas} R'(X,Y)Z &=
R(X,Y)Z+(\nabla_X\te)(Y,Z)-(\nabla_Y\te)(X,Z) + [\te(X),\te(Y)]Z\\
\label{eq:torcoes} T'(X,Y)&= T(X,Y)+\te(X)Y-\te(Y)X,
\end{align}
for each $X,Y,Z\in\bfG(TM)$. Where $R'$ and $T'$ denote the curvature and torsion tensors of $\nabla'$, respectively; $R$ and $T$ denote the curvature and torsion tensors of $\nabla$, respectively.
\end{obs}

\end{subsection}
\end{section}

\begin{section}{Infinitesimally homogeneous manifolds\/}\label{sectioninfinitesimally}


Let $\tripla$ be an $n$-dimensional affine manifold with $G$-structure $P$, the inner torsion of $P$ with respect to the connection $\nabla$ was introduced in
\cite{T}, this notion gives rise to a tensor field $\I^P$ on $M$ that measures the lack of compatibility of the connection $\nabla$ with $P$, since this
notion plays an important roll in this work, we present below its definition in a brief way.

For each $x\in M$, we denote by $G_x$ the Lie subgroup of $\GL(T_xM)$ consisting of {\em $G$-structure preserving\/} endomorphisms of $T_xM$. Clearly $G_x = \mathcal{I}_p(G)$, for all $p \in P_x$, so that $G_x$ is a Lie subgroup of $\GL(T_xM).$ We denote by $\mathfrak{g}_x\subset \glTx$ the Lie algebra of $G_x$. It is clear that $\Ad_p(\mathfrak{g})=\mathfrak{g}_x$, for all $p\in P_x$, where $\mathfrak{g}\subset \gl$ denotes the Lie algebra of $G$. Since $\Ad_p:\gl \to \glTx$ carries $\mathfrak{g}$ onto $\gx$; therefore, it induces an isomorphism:
\[\overline{\rm{Ad}}_p :\gl/\mathfrak{g} \to \glTx/\gx.\] Let $s:U\subset M \to P$ be a smooth local section of $P$, with $x\in U$ and set $s(x)=p$. If $\omega$ denotes the $\gl$-valued connection form on $\FR(TM)$ associated with $\nabla$ and $\ol{\omega}=s^*{\omega}$.
The map
\begin{equation}\label{eq:e2}\xymatrix{%
T_xM \ar@.@/_1.5pc/[rrr]_{\I ^P_x}\ar[r]^-{\overline{\omega}_x}&
\mathfrak{gl(n)} \ar[r]^-{\mathfrak q}& \mathfrak{gl(n)}/\mathfrak{g}
\ar[r]^-{\overline{{\rm{Ad}}}_{p}}& \glTx/\gx }
\end{equation}
does not depend on the choice of the local section $s$. The linear map $\I^P_x$ defined by
\eqref{eq:e2} is called the {\em inner torsion\/} of the $G$-structure $P$ at the point $x$ with respect to the connection $\nabla$.
It follows from \eqref{eq:tensorcriseforma}, that if $s:U\to P$ is a smooth local section with $x\in U$ and $\Gamma$ denotes the Christoffel
tensor of $\nabla$ with respect to $s$ then the inner torsion $\I^P_x$ is precisely the composition of the $\Gamma_x
:T_xM \to \glTx$ with the quotient map $\glTx\to \glTx/\gx$. This observation gives a simple way of computing inner torsions, (see \cite{T}).

The geometry of an affine manifold with $G$-structure $(M,\nabla,P)$ is described by three tensors of $M$: the torsion $T$ of $\nabla$, the curvature $R$ of $\nabla$
and the inner torsion $\I^P$. An important class of examples of affine manifolds with $G$-structure is defined by the property that these three tensors $T$, $R$ and $\I^P$
be {\em constant\/} when written in frames of the $G$-structure $P$. When this is the case, $(M,\nabla,P)$ is said to be {\em infinitesimally homogeneous\/}. This
statement is made more precise in the following definition.



\begin{df}\label{thm:tensoresconstantes}
An $n$-dimensional affine manifold with $G$-structure, $\tripla$
is said to be {\em infinitesimally homogeneous\/} if there exists
maps $R_0\in\Lin_3(\R^n,\R^n)$, $T_0\in\Lin_2(\R^n,\R^n)$ and a
linear map $\I_0:\R^n\to\gl/\mathfrak g$ such that: for every $x\in M$, every
$p\in P_x$ relates $T_0$ with $T_x$, $R_0$ with $R_x$ and $\overline\Ad_p\circ\I_0=\I^P_x\circ p$.
\end{df}
The maps $T_0$, $R_0$, $\I_0$ as refered above are called the {\em cha\-racteristic tensors\/} of the infinitesimally homogeneous manifold $(M,\nabla,P)$.

Clearly, the charac\-teristic tensors $T_0$, $R_0$, $\I_0$ of an infinitesimally homogeneous manifold $(M,\nabla,P)$ are invariant by the
action of the {\em structural group\/} $G$. Therefore, it follows from the $G$-invariance condition that the following relations hold:
\begin{gather}
\label{eq:R0}
R_0(u,v)=\Ad_g\cdot R_0(g^{-1}\cdot u,g^{-1}\cdot v);\\
\label{eq:T0}
T_0(u,v)=g\cdot T_0(g^{-1}\cdot u,g^{-1}\cdot v);\\
\label{eq:e6}
\Ad_g\big(\lambda(g^{-1}\cdot
u)\big)-\lambda(u)\in \g,
\end{gather}
for all $g\in G$, all $u, v\in \R^n$. Where $\lambda: \R^n\to\gl$ is an arbitrary lifting of $\I_0$.
Notice that relation \eqref{eq:e6} does not depend on $\lambda$. In fact, let $\lambda, \delta$ be liftings of $\I_0$. Write $\lambda = \delta+L$, where $L$ is a $\g$-valued linear map defined in $\R^n$. An easy computation shows that: \[\g \ni \Ad_g\big(\lambda(g^{-1}\cdot u)\big)-\lambda(u) = \Ad_g\big(\delta(g^{-1}\cdot u)\big)-\delta(u) + \underbrace{\Ad_g\big(L(g^{-1}\cdot u)\big)-L(u)}_{\in \g},\] for all $g\in G$, $u\in\R^n$.

By differentiating \eqref{eq:R0}, \eqref{eq:T0}, and \eqref{eq:e6} we obtain the following:
\begin{lema}\label{thm:g-invariantes1} Let $\lambda:
\R^n\to\gl$ be an arbitray lifting of $\I_0$. Then for all $L
\in \g$ and all $u,v \in \R^n$, the following conditions hold:
\begin{enumerate}
\item
$[L, R_0(u,v)]-R_0(L\cdot u,v)-R_0(u,L\cdot v) =0;$
\item
$L\circ T_0(u,v)-T_0(L\cdot u,v)-T_0(u,L\cdot v)=0;$
\item
$[L,\lambda(u)]-\lambda(L\cdot u)\in \g$.
\end{enumerate}
\end{lema}
\end{section}

\begin{section}{ Algebraic relation between the characteristic tensors\/}
It is a natural question to ask whether one can give a (local) {\em classification\/}
of infinitesimally homogeneous manifolds with prescribed group $G$ and prescribed characteristic tensors $T_0$, $R_0$, $\I_0$.
We solve this question in this paper by giving ne\-cessary and sufficient conditions for maps $T_0$, $R_0$, $\I_0$ to be the characteristic tensors of an infinitesimally homogeneous manifold. Our plan for developing the necessary condition is the following: we show that to give a classification of infinite\-simally homogeneous manifolds with prescribed group $G$ is equivalent to finding an infinitesimally homogeneous manifold without torsion whose structural group is $G$, and to give a classification of the $G$-invariant maps $\te_0\in \Lin_2(\R^n,\R^n)$. Once,  this is done, in order to obtain the aimed condition, it will suffice to consider the case of
symmetric connections (equivalently $T_0=0$). This is the purpose of this section, and the sufficient conditions will be developed in the following section.

\subsection{Covariant derivative for $G$-constant tensors\/}\label{thm:derivadatensorgconstante}

Let $\tripla$ be an homogeneous affine manifold with $G$-structure $P$. If $\I^P=0$, i.e., the covariant deri\-vative of $P$ is zero, it follows that every $G$-constant tensor is parallel with respect to $\nabla$. On the other hand, if $\nabla$ is not compatible with $P$, i.e., the covariant deri\-vative of $P$ is not zero, this is not true. In what follows we will show a simply way to calculate the covariant derivative for $G$-constan tensors on this case, i.e., when $\I^P\ne 0$.

Denoting by $\Vect{}$ the category whose objects are real finite-dimensional
vector spaces and whose morphisms are linear isomorphisms. Given a smooth
functor $\Cat F:\Vect{} \to\Vect{}$ and any object $V$ of $\Vect{}$, $\Cat F$
induces a Lie group homomorphism $\Cat F:\GL(V)\longrightarrow\GL\big(\Cat F(V)\big),$ whose differential at the identity is a
Lie algebra homomorphism that will be denoted by $
{\Cat f}:\gls(V)\longrightarrow\gls\big(\Cat F(V)\big).
$

Let $E$ be a vector bundle with typical fiber $E_0$ over $M$. Given a smooth
functor $\Cat F:\Vect{}\to\Vect{}$ we denote by $\Cat F(E)=\bigcup_{x\in M}\Cat F(E_x)$,
the vector bundle with typical fiber $\Cat F(E_0)$ obtained from $E$ by using $\Cat F$.

Given a smooth funtor $\Cat F :\Vect {}\to
\Vect {}$ we have the following:
\begin{lema}\label{thm:derivadaconstantes}
Let $\te$ be a smooth $G$-constant section of $\Cat
F(TM)$. Then
\begin{equation}\label{eq:derivadaGconstante}\nabla_v\te = \Cat f(L)\cdot\te_x,
\end{equation} for all $x\in M$, $v \in T_xM$, where $L\in \glTx$ is such that $\mathfrak I_x^P(v)=L+\gx$.
\end{lema}
\begin{proof}
Clearly $\te$ can be thought of as an $\FR(TM)$-valued $0$-form on $M$, which
is associated to a $0$-form $\phi:\FR(TM)\to \Cat F(\R^n)$ such that:
$\phi(p)=\Cat F(p)^{-1}(\te_x)$ for all $x\in M$, $p \in \FR(TM)$. Moreover
the covariant exterior differential $\Dd \phi$ is associated to the covariant
exterior differential $\Dd \te$ of $\te$ \cite{T}. More explicitly, we have:
\begin{equation}\label{eq:derivadaconstante1}
\dif \phi_p(\zeta) =\Dd\phi_p(\zeta)=\Cat
F(p)^{-1}(\Dd\te)_x\cdot v =\Cat F(p)^{-1}\nabla_v\te,
\end{equation}
for all $x \in M$, $p\in P_x$, $v \in T_xM$ and $\zeta$ a horizontal vector
such that $\dif \pip_p(\zeta) =v$, where $\pip:\FR(TM)\to M$ denotes the
canonical projection. To obtain the desired result, we must to
calculate $\dif \phi_p(\zeta)$.
If $X\in \gl$ is such that
$\ol{\Ad}_p(X+\g)=\mathfrak I_x^P(v)$
then
\[\zeta =(\dif \pip_p,\omega_p)^{-1}(v,X)-(\dif
\pip_p,\omega_p)^{-1}(0,X)=\underbrace{(\dif
\pip_p,\omega_p)^{-1}(v,X)}_{\in T_pP}-\dif \beta_p(1)\cdot X,\] where $\beta_p$
denotes the map given by the action of ${\rm GL}(n)$ on $p$. Since $\phi\mid_P$ is constant, we have:
\begin{equation}\label{eq:derivadaconstante2}\dif \phi_p(\zeta) = -\dif \phi_p\big(\dif
\beta_p(1)\cdot X\big)=\Cat f(X)\cdot \te_0.
\end{equation}
But \eqref{eq:derivadaGconstante} follows directly from equalities \eqref{eq:derivadaconstante1},
\eqref{eq:derivadaconstante2}.
\end{proof}

\begin{example}
Let $\Cat F:\Vect {}\to \Vect {}$ be the funtor defined by: \[\Cat F(V)
=\Lin_k(V;\Lin(V))\] for each object $V$ of $\Vect {}$.
Let $\tripla$ be an $n$-dimensional affine manifold with $G$-structure. If $\te_0\in \Lin_k\big(\R^n;\gl\big)$ is a $G$-constant
tensor, denoting by $\te_x$ the induced version of $\te_0$ on $T_xM$, by using \eqref{thm:derivadaconstantes} we have:
\[\nabla_v\te=[L,\te_x(\cdot,\dots,\cdot)]-\te_x(L\cdot,\cdot,\dots,\cdot)-\cdots -\te_x(\cdot,\cdot,\dots,L\cdot),
\] where $L \in \glTx$ is such that $\mathfrak I_x^P(v) =L+\gx$.
On the other hand, it is clear that an arbitrary lifting $\lambda: \R^n\to \mathfrak{gl}(n)$ of $\mathfrak I_0$, induces for all $X \in \R^n$,
a derivation $\mathcal{D}_{\lambda(X)}$ on the tensor algebra over the vector space $\R^n$, an easy computation shows that:  \[\big(\mathcal D_{\lambda(X)}\te_0\big)=\Cat f\big(\lambda(X)\big)\cdot \te_0\]
Therefore, if $\lambda$ is an arbitrary lifting of $\I_0$, given $x\in M$, $p \in P_x$ and $X\in \R^n$ such that $v=p(X)$ and $\Ad_p(\lambda(X)) =L$ we have: \[\Ad_p(\mathcal D_{\lambda(X)}\te_0)=
(\nabla_v\te)\circ (p,\dots,p).\]
\end{example}

\subsection{Infinitesimally homogeneous manifolds without torsion\/}
Let $\tripla$ be an $n$-dimensional affine manifold with $G$-structure and assume that $\nabla$ is a symmetric connection. Let $\te_0 \in \Lin_2(\R^n,\R^n)$ be a $G$-invariant skew-symmetric tensor. For each $x\in M$, we denote by $\te_x$ the induced version of $\te_0$ on $T_xM$. In view of remark \ref{thm:diferenciadeconexiones}, it is clear that $\nabla' = \nabla +\frac 12\te$ defines a connection on $M$ whose torsion is $\te$. We devote this section to prove the following.
\begin{lema}\label{thm:lemasintorsion}
With the same notation as above, if $(M,\nabla,P)$ is an infinitesimally homogeneous manifold then the triple
$(M,\nabla',P)$ is also infinitesimally homogeneous.
\end{lema}
\begin{proof}It is enough to prove that there exists tensors $T_0'$, $R_0'$, $\I_0'$ as in \ref{thm:tensoresconstantes}. We take $T'_0=\te_0$.
On the other hand, $\te$ can be identified with a smooth $\Lin(TM)$-valued covariant
$1$-tensor field on $M$. Let $s:U\to P$ be a smooth local section of $P$.
We denote by $\Gamma'$ and $\Gamma$, respectively,  the Christoffel tensor of $\nabla'$ and $\nabla$ with respect to $s$. Given $x\in U$, it is clear that $\Gamma'_x =\Gamma_x+\te_x$, by composing this with the canonical projection $\gls(T_xM)\to \gls(T_xM)/\gx$ we obtain:
\[\mathfrak I_x'^P = \mathfrak I_x^P + \mathfrak q\circ \te_x.\] Therefore, we can take $\mathfrak I'_0 = \mathfrak I_0 + \mathfrak q\circ \te_0.$
On the other hand, we denote by $R'$ and $R$, respectively, the curvature tensor of $\nabla'$ and $\nabla$. Let  $\lambda$ be an arbitrary lifting of $\I_0$, $x\in U$ and set $s(x)=p$. From \eqref{eq:curvaturas} and by using
lemma ~\ref{thm:derivadaconstantes}
we have that the following holds:
\beq R'_x(p\cdot,p\cdot)&=&R_x(p\cdot,p\cdot
)+(\Dd\te)_x(p\cdot,p\cdot)+[\te_x(p\cdot),\te_x(p\cdot)] \nonumber \\
&=&\Ad_p\circ\big( R_0(\cdot,\cdot) +\Alt\big(\mathcal
D_{\lambda(\cdot)}\te_0\big)\cdot+[\te_0(\cdot),\te_0(\cdot)]\big).\nonumber\eeq
Therefore, in order to obtain the desired result we can take \[R'_0 = R_0 + \mathcal D
\te_0+[\te_0,\te_0].\]
 \end{proof}

\subsection{The necessary conditions\/}\label{relationalgebric}
We are now ready to give necessary conditions which must be satified by the characteristic tensors of an infinitesimally homogeneous manifold. To
do this, throughout the subsection we consider a fixed $n$-dimensional
infinitesimally homogeneous manifold $\tripla$ with structural group $G$. From
lemma \ref{thm:lemasintorsion} it fo\-llows that we may assume without loss of
genera\-lity that $\nabla$ is a symmetric connection with curvature $R$. We
denote by $R_0, \I_0$  the characteristic tensor of $\tripla$. Clearly, a
necessary condition is that $R_0, \I_0$ are $G$-invariant.

Let $\omega$ be the $\mathfrak{gl}(n)$-valued connection form on $\FR(TM)$
associated with $\nabla$, let $\Omega$ be its curvature form and let
$\theta$ be the canonical form of $\FR (TM)$. Given a smooth local frame
$s:U\to P$ then, setting $\ol \omega =s^*(\omega)$, $\ol{\Omega}=s^*\Omega$, $\ol{\theta}=s^*\theta$, we have: \[\ol \Omega = \dif \ol \omega + \ol \omega \wedge \ol \omega,\;\;\;\dif \ol{\theta}=-\ol{\omega}\wedge \ol{\theta}.\]
Moreover, the infinitesimal homogenity implies that: \[\ol \Omega_x(X,Y)= s(x)\circ R_x(X,Y) \circ s(x)^{-1}=R_0\big(s(x)^{-1}X,s(x)^{-1}Y\big),\]
\[ \mathfrak q \circ \ol{\omega}_x = \ol{\Ad}_{s(x)^{-1}}\circ \I_x^P=\I_0 \circ \ol{\theta},\] for all $x\in U$, $X, Y \in T_xM$, where $\mathfrak{q}:\mathfrak{gl}(n)\to \mathfrak{gl}(n)/\mathfrak{g}$ denotes the canonical projection and $\mathfrak{g}$ denotes the Lie algebra of $G$.

Clearly when the linear map $\I^P$ vanishes, $\ol \Omega$ is a $\mathfrak{g}$-valued $2$-form on $M$. Under the previous conditios, in order to handle the
general case in which $P$ is not compatible with $\nabla$ we get:
\beq \mathfrak q\circ \ol{\Omega}
&=& \dif (\mathfrak q \circ \ol{\omega})+
\mathfrak q \circ \ol{\omega}\wedge \ol{\omega} \nonumber \\
&=& \dif (\I_0 \circ \ol{\theta})+ \mathfrak q \circ \ol{\omega}\wedge \ol{\omega} \nonumber \\
&=& \I_0\circ \dif \ol{\theta} +\mathfrak q \circ
\ol{\omega} \wedge \ol{\omega} \nonumber\\
&=& \label{eq:proyecciondelaforma}-\I_0 \circ (\ol{\omega} \wedge \ol{\theta}) + \mathfrak q \circ\ol{\omega} \wedge \ol{\omega}.
\eeq
Given $x\in U$, let $\wt{\Gamma}: \R^n \to \gl$ be the map defined by requiring
the diagram
\[
\xymatrix{T_xM \ar[r]^-{\Gamma_x}  \ar[rd]_{\ol{\omega}_x}& \mathfrak{gl}(T_xM)\\
\R^n \ar[r]_-{\wt{\Gamma}} \ar[u]^{s(x)}  & \mathfrak{gl}(n) \ar[u]_{\Ad_{s(x)}} }
\]
to be commutative. Therefore,  $\I_0 =\mathfrak q \circ \wt{\Gamma}$ and substituting
in \eqref{eq:proyecciondelaforma} we obtain the following relation:
\begin{equation*} \ol{\Omega}_x+\wt{\Gamma} \circ
(\ol{\omega}_x\wedge\ol{\theta}_x)-\ol{\omega}_x\wedge
\ol{\omega}_x \in \g.\end{equation*}

Thus, given vectors $u, v \in \R^n$ the relation above can be written as:
\be\label{eq:relentreReI}
R_0(u,v)-[\wt{\Gamma}(u),\wt{\Gamma}(v)]+\wt{\Gamma}\big(\wt{\Gamma}(u)v-\wt{\Gamma}(v)u\big)
\in \g.\ee
This relation does not depend on the choice of $\tilde{\Gamma}$. Namely,
let $\lambda$ be an arbitrary lifting of $\I_0$ and $\delta$ be a $\g$-valued linear map in $\R^n$
such that $\wt{\Gamma}= \lambda+\delta$. By replacing this into \ref{eq:relentreReI}, we obtain
\begin{equation}\label{eq:relationRandI} \g \ni
R_0(u,v)-[\lambda(u),\lambda(v)]+\lambda\big(\lambda(u)v-\lambda(v)u\big)
+ \mathcal A(\delta) + \mathcal B(\delta),\end{equation} where

\begin{align*}
\mathcal
A(\delta)&=\big([\delta(v),\lambda(u)]-\lambda(\delta(v)\cdot u)\big)- \big([\delta(u),\lambda(v)]-\lambda(\delta(u)\cdot v)\big),\\
\mathcal B(\delta)&=\delta\big(\wt{\Gamma}(u)v - \wt{\Gamma}(v)u\big)-[\delta(u),\delta(v)].\\
\end{align*}

So that Lemma~\ref{thm:g-invariantes1} guarantees that $\mathcal A(\delta)
\in \g$; moreover, $\mathcal B(\delta) \in \g$ because $\delta$ is a $\g$-valued linear map. Therefore for an arbitrary lifting $\lambda$ of $\I_0$ the following relation
holds:
\[R_0(u,v)-[\lambda(u),\lambda(v)]+\lambda\big(\lambda(u)v-\lambda(v)u\big)\in \g,\] this shows the independence on the lifting; hence we have proved the following:

\begin{teo}\label{thm:relentreReI}
Let $M$ be an $n$--dimensional differentiable manifold, $G$ a Lie subgroup of
$\GL(n)$ with Lie algebra $\mathfrak{g}$ and assume that $M$ is endowed with a symme\-tric connection $\nabla$ and a $G$--structure
$P\subset \FR(T M )$. Assume that $\tripla$ is an infinitesimally homogeneous manifold with characteristic tensors $R_0$, $\I_0$. Then given an arbitrary lifting $\lambda$ of $\I_0$, the following relation holds:
\[R_0(u,v)-[\lambda(u),\lambda(v)]+\lambda\big(\lambda(u)v-\lambda(v)u\big)\in
\g,\] for all $u,v \in \R^n$.
\end{teo}
\end{section}

\begin{section}{Infinitesimally homogeneous manifolds with prescribed group and prescribed characteristic tensors}\label{inverseproblem}

We devote this section to obtain sufficient conditions for maps $T_0$, $R_0$, $\I_0$ to be the characteristic tensors of an infinitesimally homogeneous manifold. 
Therefore, in this section we will consider fixed a real finite-dimensional vector
space $\m$, a Lie subgroup $H \subset \GL(\m)$ with Lie algebra $\h \subset \gls(\m)$ and $H$--invariant maps $R_0 \in \Lin_2\big(\m,\gls(\m)\big)$, $\I_0:\m\to \gls(\m)/\h$. As we said above, our goal is to obtain conditions for the maps $R_0$, $\I_0$ to be the characteristic tensors of an infinitesimally homogeneous manifold $(M,\nabla,P)$.

Let $\lambda:\m\to \gls(\m)$ be an arbitrary lifting of $\I_0$. As in Section \ref{sectioninfinitesimally},  by using the $H$--invariance of $\I_0$ we conclude that the following relation holds:
\be
\label{eq:inverso2}[L,\lambda(X)]-\lambda(L\cdot X)\in \h,
\ee for all $L\in \h$, all $X, Y \in \m$.
An analogous relation to \eqref{eq:relentreReI} is:
\be\label{eq:inverso3}R_0(X,Y)-[\lambda(X),\lambda(Y)]+\lambda\big(\lambda(X)Y-\lambda(Y)X\big) \in \h\ee for all $X,Y \in \m$.
Neither relation \eqref{eq:inverso2} nor relation \eqref{eq:inverso3} do not depend on the choice of $\lambda$.

Assuming that \eqref{eq:inverso3} holds, we have the following:
\begin{df}\label{thm:colchete}
Setting $\A = \h \oplus \m$. We endow $\A$ with a bracket operation which is defined below. For each $X,Y \in \m$, each $L, T \in \h$ we set:
\begin{enumerate}
\item $[X,Y]^{\m} = \lambda(X)\cdot Y-\lambda(Y)\cdot X$;

\item
$[X,Y]^{\h} =R_0(X,Y)+\lambda\big(\lambda(X)\cdot Y-\lambda(Y)\cdot X\big)- [\lambda(X),\lambda(Y)]$;

\item
$[L,X]^{\m}=L\cdot
X$;

\item
$[L,X]^{\h}=[L,\lambda(X)]-\lambda(L\cdot X)$;

\item
$[L,T]$ is the Lie bracket of $\h$;

\item
$[L,X]=-[X,L]$.
 \end{enumerate}
\end{df}

We will prove that the vector space $\A$ endowed with the bracket operation as above is a Lie algebra. Before we proceed, we will present some algebraic preli\-minaries.



\begin{df}
We say that the map $R_0$ satisfies the {\em Bianchi identities} if the fo\-llowing equalities hold:
\begin{itemize}
\item[$(B_1)$]
$\mathfrak{S}R_0(X,Y)\cdot Z =0$;
\item[$(B_2)$]
$\mathfrak{S}\big(\mathcal D_{\lambda(X)}R_0\big)(Y,Z) =0$.
\end{itemize}
Where for $ X\in \m$, $\mathcal{D}_{\lambda(X)}$  denotes the derivation on the tensor algebra over the vector space $\m$ induced by $\lambda(X)$ and $\mathfrak{S}$ denotes the sum over all cyclic permutations of $X, Y, Z.$
\end{df}


\begin{obs}\label{rem:4}
For $X, Y, Z \in \m$ and $L\in \h$ we will use the next notation:
\begin{align*}
\mathcal S_{[L,X,Y]}&=[L,\lambda(X)]\cdot Y -\lambda(Y)\cdot
(L\cdot X). \\
\mathcal T_{[X,Y,Z]}&=[\lambda(X),\lambda(Y)]\cdot
Z-\lambda(Z)\cdot[X,Y]^{\m}.
\end{align*}
Thus, it is not difficult to see that:
\beq \mathcal S_{[L,X,Y]}-\mathcal
S_{[L,Y,X]}
\label{eq:exprssaoS}&=& L\big([X,Y]^{\m}\big). \eeq
We can also easily see that: \be\label{eq:expressaoT} \mathfrak{S}
\mathcal T_{[X,Y,Z]}=0.
\ee
\end{obs}

\begin{obs}
For $X, Y, Z \in \m$ by using the Bianchi identities we obtain:
\be\label{eq:sumaparaR}
\mathfrak{S}\Big([\lambda(Z),R_0(X,Y)]-R_0\big([X,Y]^{\m},Z\big)\Big)=0.  \ee
\end{obs}
\newpage
\begin{lema}\label{thm:algebradeLie}
Using the same notations and terminology as above, suppose that the $H$--invariant maps $R_0, \I_0$ satisfy the following conditions
\begin{enumerate}
\item
$R_0$ is skew-symmetric;
\item
given an arbitrary lifting $\lambda:\m\to \gls(\m)$ of $\I_0$, the map $R_0$ satisfies the Bianchi identities and the relation \eqref{eq:inverso3} holds.
\end{enumerate} Then the vector space $\A= \h\oplus \m$ endowed with the bracket operation $[\cdot,\cdot]$, defined as in \eqref{thm:colchete}, is a Lie algebra.
\end{lema}

\begin{proof}[Proof of Lemma \ref{thm:algebradeLie}]
Since $[\cdot,\cdot]$ is skew-symmetric, it is enough to show that satifies  the Jacobi identity. To do that,  we divide the proof in three cases. First we consider the case that $L,T \in \h$, $X\in \m$. It  follows from definition~\ref{thm:colchete}
that:

\be\label{eq:inverso8}\big[[X,L],T\big] =
-\big[[L,\lambda(X)],T\big] -\lambda\big(T(L\cdot X)\big)+T(L\cdot
X).
\ee

Interchanging $T$ and $L$ in \eqref{eq:inverso8} we get:

\be\label{eq:inverso9} \big[[T,X],L\big]=\big[[T,\lambda(X)],L\big]
+\lambda\big(L(T\cdot X)\big)-L(T\cdot X).
\ee
On the other hand, it follows from definition \ref{thm:colchete} that:

\be\label{eq:inverso12} \big[[L,T],X\big] =
\big[[L,T],\lambda(X)\big]-\lambda\big([L,T]\cdot X\big) +
[L,T]\cdot X.
\ee
The conclusion follows from \eqref{eq:inverso8}, \eqref{eq:inverso9} and \eqref{eq:inverso12} by applying the Jacobi identity in $\gls(\m)$.

Now we consider the case that $X, Y \in \m$, $L\in \h$. In this case, we get:

\begin{align}
\label{eq:casom1}\big[[X,Y],L\big]^{\m}&=-L\big([X,Y]^{\m}\big)\\
\label{eq:casoh1}\big[[X,Y],L\big]^{\h}&=\big[[\lambda(X),\lambda(Y)],L\big]+\lambda\big(L\cdot[X,Y]^{\m}\big)-[R_0(X,Y),L],
\end{align}  and using remark~\ref{rem:4} we obtain:

\beq
\label{eq:casom2}  \big[[Y,L],X]^{\m} &=& -\mathcal S_{[L,Y,X]}\\
\label{eq:casoh2} \big[[Y,L],X\big]^{\h} &=&\big[[\lambda(Y),L],\lambda(X)\big]
+\lambda\big(\mathcal S_{[L,Y,X]}\big)-R_0(X,L\cdot Y).
\eeq Interchanging $X$ and $Y$ in \eqref{eq:casom2}, \eqref{eq:casoh2} we get:

\beq
\label{eq:casom3}\big[[L,X],Y]^{\m}&=& \mathcal S_{[L,X,Y]}.\\
\label{eq:casoh3}\big[[L,X],Y\big]^{\h}&=&\big[[L,\lambda(X)],\lambda(Y)\big]
-\lambda\big(\mathcal S_{[L,X,Y]}\big)+R_0(Y,L\cdot X).
\eeq  It follows from \eqref{eq:casom1}, \eqref{eq:casom2} and \eqref{eq:casom3} by
using \eqref{eq:exprssaoS} that:

\[\mathfrak{S}\big[[X,Y],L\big]^{\m}=0.\] On the other hand, it follows from \eqref{eq:casoh1}, \eqref{eq:casoh2} and \eqref{eq:casoh3} by using \eqref{eq:exprssaoS}, \eqref{eq:sumaparaR}
and the Jacobi identity in $\gls(\m)$ that:
\[\mathfrak{S}\big[[X,Y],L\big]^{\h}=0.\]

Finally, we consider the case $X, Y, Z \in \m$. It follows directly from definition \ref{thm:colchete} that:

\[\mathfrak{S}\big[[X,Y],Z\big]^{\m}=0\]
For the $\h$ component we have:

\beq
\big[[X,Y],Z\big]^{\h}&=&\big[[\lambda(X),\lambda(Y)],\lambda(Z)\big]-R_0\big([X,Y]^{\m},Z\big)\nonumber\\&-&[R_0(X,Y),\lambda(Z)]-\lambda\big(\mathcal
T_{[X,Y,Z]}-R_0(X,Y)Z \big).\nonumber
\eeq
It follows from \eqref{eq:expressaoT} and \eqref{eq:sumaparaR}  by using the Jacobi identity in $\gls(\m)$ that:\[\mathfrak{S}\big[[X,Y],Z\big]^{\h}=0.\]
\end{proof}

We now must prove that the Lie bracket defined in \ref{thm:colchete} does not depend on the choice of $\lambda$. In fact, if $[\cdot,\cdot]_{\lambda}$ denotes the Lie Bracket in $\A$ obtained by using
the arbitrary lifting $\lambda$ of $\I_0$, given another lifting $\tilde{\lambda}$ there exists a linear map $\delta:\m\to\h$ such that $\lambda = \tilde{\lambda}+\delta$. The map $\varphi:\A\to\big(\A,[\cdot,\cdot]_{\tilde{\lambda}}\big)$ defined by the matrix:
\[\left(\begin{array}{cc} \Id_{\h}&\delta\\
0&\Id_{\m}\\
\end{array}\right),\] is an isomorphism of vector spaces, moreover, a direct computation shows that $[\cdot,\cdot]_{\lambda}=\varphi^*[\cdot,\cdot]_{\tilde{\lambda}}$ so that $\varphi$ is an isomorphism of Lie algebras. Which shows the assertion.

\subsection{Existence of an infinitesimally homogeneous manifold}
The main goal of this subsection is to show the existence of an infinitesimally homogeneous manifold with prescribed structural group and prescribed characteristic tensors. To do this, let $\m$ be a real finite-dimensional vector space, let $H\subset \GL(\m)$ be a Lie subgroup with Lie algebra $\h\subset \gls(\m)$. Let $R_0 \in \Lin_2\big(\m,\gls(\m)\big)$, $\I_0:\m\to \gls(\m)/\h$, be maps satisfying the following conditions:
\begin{enumerate}
\item
$R_0, \I_0$ are $H$-- invariants;
\item
$R_0$ is skew--symmetric;
\item
given an arbitrary lifting $\lambda:\m\to \gls(\m)$ of $\I_0$, $R_0$ satisfies the Bianchi identities and the relation \eqref{eq:inverso3} holds.
\end{enumerate}  Now we are going to obtain an infinitesimally homogeneous manifold with structural group $H$ whose characteristic tensor are $R_0$, $\I_0$.
It follows from Lemma \ref{thm:algebradeLie} that the vector space $\A=\h\oplus \m$ endowed with the bracket defined on \ref{thm:colchete} is a Lie algebra.

\noindent Let $\ol{\lambda}:\A=\h\oplus \m\to \gls(\m)$ be a map defined by:
\begin{equation}\label{eq:lambdabarra}
\ol{\lambda}(X)=\begin{cases}
\lambda(X),&\text{se $X\in\m$},\\
\ol{\ad}_X,&\text{se $X\in\h$},
\end{cases}
\end{equation} for each $X\in \A$, where
$\ol{\ad}$ denotes the isotropic representation of $\h$ on $\m$, more
precisely $\ol{\ad}_X(Y) = \pro \m\big([X,Y]\big)=X(Y)$ for all
$X\in \h$, $Y\in \m$.

\begin{lema}\label{thm:propriedadelambda}
If $L\in \h$ and $\mathfrak X\in \A$. Then \[\big[\ol{\lambda}(L),\ol{\lambda}(\mathfrak
X)\big]=\ol{\lambda}\big([L,\mathfrak X]\big).\]
\end{lema}
\begin{proof}
We set $\mathfrak X = T+X$, for $T \in \mathfrak h$, $X\in \mathfrak m$.
\beq \ol{\lambda}\big([L,\mathfrak X]\big)
&=& \ol{\ad}_{[L,T]}+\ol{\ad}_{\pro \h\big([L,X]\big)}+
\lambda(L\cdot X)\nonumber \\
&=& [\ol{\ad}_L,\ol{\ad}_T]+[\ol{\ad}_L,\lambda(X)]\nonumber \\
&=& [\ol{\lambda}(L),\ol{\lambda}(\mathfrak X)].\nonumber
\eeq
\end{proof}

Let $A$ be a Lie group such that $T_1A = \A$. Let $M'\subset A$ be a submanifold of $A$ throught $1$ such that $T_1M'=\m$.
Let $\pro \m^L$ be the left invariant $1$-form on $A$ induced by the linear projection $\pro \m:\A = \h\oplus \m \to \m$. Setting $\ol{\kappa}=\pro \m^L|_{M'}$ then:  \[\ol{\kappa}_1(X)=\pro \m^L(X)= \pro \m(X) =X\] for all $X\in \m$. Let $M$ be a neighborhood of $1$ in $M'$ such that for all $x\in M$ the map $\ol{\kappa}_x:T_xM \to \m$ is a linear isomorphism. Then, the map $s:M\to\FR_{\m}(TM)$ defined by $s(x) = \ol{\kappa}_x^{-1}:\m\to T_xM$, for all $x\in M$ gives us a global section of the $\GL(\m)$-principal bundle $\FR_{\m}(TM)$ over $M$. Given $x\in M$, the
set \[P_x =s(x)\cdot H=\{s(x)\circ h:h \in H\},\] is an $H$-structure on $T_xM$ and $\displaystyle P=\bigcup_{x\in M}P_x$ defines an $H$-structure on $M$.

\noindent To construct $\nabla$, let $\ol{\lambda}^L$ the left invariant $1$-form on $A$ induced by the linear map $\ol{\lambda}$ defined in \eqref{eq:lambdabarra}. Setting $\ol{\omega}=\ol{\lambda}^L|_M$, it is clear
that $\ol{\omega}$ is a $\gls(\m)$-valued smooth $1$-form on $M$. Let  $\omega$ be the unique $\gls(\m)$-valued $1$-form on $\FR_{\m}(TM)$ such that $s^*\omega=\ol{\omega}$. Then $\omega$ is a connection form on $\FR_{\m}(TM)$, (see \cite{T}).
\newline

So far, we have obtained an affine manifold with $H$-structure $\tripla$, where $\nabla$ denotes the linear connection associated with the connection form $\omega$.
We claim that $\tripla$ is an infinitesimally homogeneous manifold whose characteristic tensors are $R_0$, $\I_0$. In fact,
given $x\in M$ and $X\in T_xM$, we have:
\[\ol{\omega}_x(X)=\ol{\lambda}^L_x(X)=\ol{\lambda}(x^{-1}\cdot
X)=\underbrace{\ol{\ad}_{\pro \h(x^{-1}\cdot X)}}_{\in
\h}+\lambda\big(\pro \m(x^{-1}\cdot X)\big),\]
therefore, in the quotient $\mathfrak{gl}(\m)/\h$ the following equality holds:
\[\ol{\omega}_x(X)=\lambda\big(\pro \m(x^{-1}\cdot X)\big);\] clearly
$\pro \m(x^{-1}\cdot X) = \ol{\kappa}_x(X)=s(x)^{-1}\cdot X$. Thus we have:
\[\I_x^P (X) =\ol{\Ad}_{s(x)}\big(\mathfrak q\circ \lambda \circ
s(x)^{-1}\cdot X\big)=\ol{\Ad}_{s(x)}\big(\I_0\circ s(x)^{-1}\cdot
X\big).\]
On the other hand, we set $\ol{\Omega}=s^*\Omega$, where $\Omega$ denotes the curvature form of $\omega$. For each  $x\in M$, $X, Y\in T_xM$.
Setting $x^{-1}\cdot X
=L+\ol{\kappa}_x\cdot X$, $x^{-1}\cdot Y =T+\ol{\kappa}_x\cdot Y$, for $L,T \in \h$. It follows from Lemma~\ref{thm:propriedadelambda} that:
\beq
-\ol{\omega}_x\big([X,Y]\big)&=&-\ol{\lambda}\big([L,T+\ol{\kappa}_x\cdot
Y]+[\ol{\kappa}_x\cdot X,T]+[\ol{\kappa}_x\cdot X,\ol{\kappa}_x\cdot
Y]\big)\nonumber \\
&=& -\big[\ol{\lambda}(L),\ol{\lambda}(T+\ol{\kappa}_x\cdot
Y)\big]-\big[\ol{\lambda}(\ol{\kappa}_x\cdot
X),\ol{\lambda}(T)\big]\nonumber\\&-&\ol{\lambda}\big[\ol{\kappa}_x\cdot
X,\ol{\kappa}_x\cdot Y\big]. \nonumber\eeq
Moreover:
\beq
\big[\ol{\omega}_x(X),\ol{\omega}_x(Y)\big]
&=&\big[\ol{\lambda}(L),\ol{\lambda}(T+\ol{\kappa}_x\cdot
Y)\big]+\big[\ol{\lambda}(\ol{\kappa}_x\cdot
X),\ol{\lambda}(T)\big]\nonumber\\&+&\big[\ol{\lambda}(\ol{\kappa}_x\cdot
X),\ol{\lambda}(\ol{\kappa}_x\cdot Y)\big].\nonumber \eeq Since
\[\ol{\Omega}_x(X,Y)=\dif \ol{\omega}_x(X,Y)
+\big[\ol{\omega}_x(X),\ol{\omega}_x(Y)\big]
=-\ol{\omega}_x\big([X,Y]\big)+\big[\ol{\omega}_x(X),\ol{\omega}_x(Y)\big],\]
it follows from the previous equalities that:
\beq
\ol{\Omega}_x(X,Y)&=&-\ol{\lambda}\big[\ol{\kappa}_x\cdot
X,\ol{\kappa}_x\cdot Y\big]+\big[\ol{\lambda}(\ol{\kappa}_x\cdot
X),\ol{\lambda}(\ol{\kappa}_x\cdot Y)\big]\nonumber \\
&=& R_0(\ol{\kappa}_x\cdot
X,\ol{\kappa}_x\cdot Y).\nonumber
\eeq
Which shows the claim. The following Theorem summarizes all subsection:

\begin{teo}
Let $\m$ be a real finite-dimensional vector space, let $H\subset \GL(\m)$ be a Lie subgroup with Lie algebra $\h\subset \gls(\m)$. Let $R_0 \in \Lin_2\big(\m,\gls(\m)\big)$, $\I_0:\m\to \gls(\m)/\h$, be maps satisfying the following conditions:
\begin{enumerate}
\item
$R_0, \I_0$ are $H$-- invariants;
\item
$R_0$ is skew--symmetric;
\item
given an arbitrary lifting $\lambda:\m\to \gls(\m)$ of $\I_0$, the map $R_0$ satisfies the Bianchi identities and the relation
\[R_0(X,Y)-[\lambda(X),\lambda(Y)]+\lambda\big(\lambda(X)Y-\lambda(Y)X\big) \in \h\] holds.
\end{enumerate}
Then there exists an infinitesimally homogeneous manifold $\tripla$ with structural group $H$, whose cha\-racteristic tensors are $R_0$, $\I_0$.
\end{teo}
\end{section}

\end{document}